\input amstex 
\documentstyle{amsppt}
\input bull-ppt
\keyedby{bull350/kmt}

\topmatter
\cvol{28}
\cvolyear{1993}
\cmonth{January}
\cyear{1993}
\cvolno{1}
\cpgs{109-115}
\title Singularities of the Radon transform \endtitle
\author A. G. Ramm and A. I. Zaslavsky\endauthor
\shorttitle{Singularities of the Radon transform}
\address Mathematics Department,
Kansas State University, Manhattan, Kansas 
66506-2602\endaddress
\ml ramm\@ksuvm.ksu.edu\endml
\address Department of Mathematics, Technion-Israel 
Institute of Technology,
32000 Haifa, Israel\endaddress
\ml mar9315\@technion.technion.ac.il\endml

\date March 11, 1992\enddate
\abstract Singularities of the Radon transform of a 
piecewise smooth function
$f(x)$, $x\in R^n$, $n\geq 2$, are calculated. If the 
singularities of the
Radon transform are known, then the equations of the 
surfaces of discontinuity
of $f(x)$ are calculated by applying the Legendre 
transform to the functions,
which appear in the equations of the discontinuity 
surfaces of the Radon
transform of $f(x)$; examples are given. Numerical aspects 
of the problem of
finding discontinuities of $f(x)$, given the 
discontinuities of its Radon
transform, are discussed.\endabstract
\endtopmatter

\document
\heading I. Introduction\endheading
\par Let $f(x)$ be a compactly supported function, $D$ be 
its support, and
$\Gamma=\partial D$ be a union of finitely many $C^\infty$ 
hypersurfaces
$\Gamma_1,\dots,\Gamma_s$ in general position, each of 
which can be written in
local coordinates as
$$x_n=g(x'),\qquad x'=(x_1,\dots,x_{n-1}),n\geq 2,$$
where $g(x')\in C^\infty$, $f(x)\in C^\infty(D)$, 
$f(x)|_\Gamma\geq c>0$. The
discontinuity surface of $f(x)$ is $\Gamma$, the boundary 
of $D$. We assume
that the rank of the Hessian 
$g_{ij}(x)\coloneq\partial^2g/\partial x_i\partial
x_j$ is constant on each of $\Gamma_j$, $1\leq j\leq s$.
\par Define the Radon transform (RT) of $f(x)$ by the 
usual formula [GGV] $\hat
f(p,\alpha)=\int_{\Bbb R^n}f(x)\delta(p-\alpha\cdot 
x)\,dx$, where $\delta$ is
the delta-function. It is well known that $\hat f(\lambda
p,\lambda\alpha)=|\lambda|^{-1}\hat f(p,\alpha)$, 
$\lambda\in R^1$,
$\lambda\not=0$. Consider the integral
$$R(p,\alpha;f)\coloneq\int_{l_{\alpha p}}f(x)\mu(dx),\tag 
1$$
where $l_{\alpha p}$ is the plane $\alpha\cdot x-p=0$, 
$\alpha\in\Bbb R^n$,
$p\in\Bbb R^1$, and $\mu(dx)$ is the Lebesgue measure on 
$l_{\alpha p}$. One
has $R(p,\alpha;f)=\hat f(p/|\alpha|,\alpha^0)$,
$\alpha^0\coloneq\alpha|\alpha|^{-1}$, so that 
$R(p,\alpha;f)=|\alpha|\hat
f(p,\alpha)$, $|\alpha|=(\alpha^2_1+\cdots+
\alpha^2_n)^{1/2}$.
\par The problems we are interested in are: (P1) Find the 
singularities
of
$R(p,\alpha;f)$; and (P2) Find the surface $\Gamma$ of 
discontinuity of $f(x)$
given the singularities of $R(p,\alpha;f)$.
\par No results concerning (P2) were known. In [N] one can 
find an estimate of
the norm of $\hat f(p,\alpha)$ in Sobolev spaces. This 
result does not give
information about (P1) and (P2). In [P] there is a result 
given without proof,
which has a relation to (P1). Our result is more general. 
In $[Q]$ it is
mentioned that the values $(\alpha\:p)$, such that 
$l_{\alpha p}$ is tangent to
$\Gamma$, play a special role. This observation is made 
quantitative in our
Theorem 1.  Our results are useful for inversion of 
incomplete tomographic data
$[R2]$.
\par The basic results are formulated in \S II. They give 
solutions of the
problems (P1) and (P2). Actually, more general problems 
are solved;
particularly, finite smoothness of $f(x)$ and $\Gamma$ is 
allowed, the role of
the intersections of $\Gamma_j$ in the study of the 
singularities of
$R(p,\alpha;f)$ is clarified, etc. In \S III proofs are 
sketched. In \S IV
examples are given. In \S V numerical aspects of problem 
(P2) are discussed.
\par We conclude this introduction by an outline of our 
ideas. First, we
describe the behavior of $R(p,\alpha;f)$ in a neighborhood 
of the set $Q_f$
which is the set of singularities of $R(p,\alpha;f)$. 
Second, we prove that, in
general, there is an equation of the set $Q_f$ which is of 
the form
$q=h(\beta)$, $\beta\in R^{n-1}$, so that $Q_f$ is a 
hypersurface. Third, we
prove that the function $g(x')$ (in the equation of 
$\Gamma)$ is the Legendre
transform of the function $h(\beta)$ (in the equation of 
$Q_f)$. Fourth, we
describe some geometric properties of $Q_f$.
\par Our results give a theoretical basis for the solution 
of the practically
important problem in nondestructive evaluation and remote 
sensing, the problem
of finding the discontinuities of a function from the 
knowledge of its RT.
\heading II. Formulation of the results\endheading
\par The RT, defined by formula (1), is a function on the 
projective space
$\Bbb R\Bbb P_n$, and we take $R(1,0;f)\coloneq 0$ for 
compactly supported $f$.
Let $Q_f$ denote the set of the points $(\alpha\:p)$ in 
this projective space,
which correspond to the planes $l_{\alpha p}$ tangent to 
$\Gamma=\partial D$.
We say that $l_{\alpha p}$ is tangent to $\Gamma$ at a 
point $x\in
B_m\coloneq\bigcap^m_{j=1}\Gamma_j$, if $l_{\alpha p}$ is 
not transversal to
$B_m$ at the point $x$.
\bh {\rm1}\endbh
Our first result is the following theorem in which the 
description of the
singularities of $R(p,\alpha;f)$ is given.
\par Let $l_{\alpha p}$ be tangent to $\Gamma$ at the 
point $\overline x$. We
claim that if $\overline\alpha$ is generic, then the set 
$Q_f$ is a smooth
hypersurface in a neighborhood $U$ of 
$(\overline\alpha\:\overline p)$. If $A$
is a symmetric matrix with real-valued entries, then its 
inertia index
(inerdex) is defined to be the number of its negative 
eigenvalues. Consider first
the case when $\Gamma$ consists of one surface.
\thm{Theorem 1} There exists an equation 
$\zeta(\alpha\:p)=0$,
$\nabla\zeta\not=0$ in $U$, which defines $Q_f$ in $U$, 
and two $C^\infty$
functions $r_1$ and $r_2$ in $U$ such that
$$R(p,\alpha;f)=\cases\zeta^{(n-1)/2}_+r_1+
r_2,&\text{if\quad {\rm In} is even},\\
\zeta^{(n-1)/2}(\operatorname{ln}|\zeta|)r_1+
r_2,&\text{if\quad {\rm In} is
odd}.\endcases\tag 2$$
Here $I$ is the inerdex of the matrix $z_{kj}$, where 
$z_{kj}$ is the Hessian
of the function $z=(\overline\alpha\cdot x-\overline 
p)/|\overline\alpha|$ on $\Gamma$
at the point $\overline x$ and $z_+=\max(z,0)$.
\ethm
\par If $\overline x\in B_m$ and $(\overline 
\alpha\:\overline p)$ is generic,
then the following result holds.
\thm{Theorem $\boldone'$} There exists $\zeta(\alpha\:p)$, 
$\nabla\zeta\not=0$ in $U$,
such that the equation $\zeta(\alpha\:p)=0$ is the 
equation of $Q_f$ in $U$,
and two $C^\infty$ functions $r_1$ and $r_2$ in $U$, such 
that
$$R(p,\alpha;f)=\cases\zeta^{(n+m-2)/2}_+r_1+r_2,&\text{if 
}\ I(n+m-1)\text{ is
even},\\
\zeta^{(n+m-2)/2}(\operatorname{ln}|\zeta|)r_1+
r_2,&\text{if }\ I(n+m-1)\text{ is odd}.\endcases
\tag"$(2')$"$$
\ethm
\par In [RZ1] the constant $r_1(\overline 
\alpha\:\overline p)$ is calculated.
In [RZ2] this result is used for a derivation of the 
asymptotics of the Fourier
transform of a piecewise smooth function.
\bh {\rm2}\endbh
Let us define the Legendre transform of a function $g(y)$, 
$y\in R^{n-1}$ in a
neighborhood $U_{\overline y}$ of a point $\overline y$ at 
which the matrix
$g_{ij}(y)\coloneq\partial^2g/\partial y_i\partial y_j$ is 
nondegenerate, i.e.,
$\operatorname{det}g_{ij}(y)\not=0$ in $U_{\overline y}$. 
Define $Lg\coloneq
h(\beta)\coloneq\beta\cdot y-g(y)$, where the dot stands 
for the inner product
and $y=y(\beta)$ is the unique solution of the equation 
$\beta=\nabla g(y)$ in
a neighborhood $U_{\overline \beta}$ of the point 
$\overline \beta=\nabla
g(\overline y)$. One can prove that if $g\in 
C^l(U_{\overline y})$, $l\geq 2$,
and $\operatorname{det}g_{ij}(y)\not=0$ in $U_{\overline 
y}$, then $h(\beta)\in
C^l(U_{\overline \beta})$.
\par It is known that under our assumptions $Lh=g(y)$, 
i.e., the Legendre
transform is involutive: $g(y)=\beta\cdot y-h(\beta)$, 
where $\beta=\beta(y)$
is the unique solution to the equation $y=\nabla 
h(\beta)$, $\beta\in
U_{\overline\beta}$. One can prove that 
$\operatorname{det}h_{ij}(\beta)\not=0$
in $U_{\overline\beta}$ if 
$\operatorname{det}g_{ij}(y)\not=0$ in $U_{\overline
y}$; moreover, the matrix $h_{ij}(\beta)$ is inverse to 
$g_{ij}(y)$, where
$\beta=\beta(y)$. Recall that $\Gamma$ is a union of 
hypersurfaces $\Gamma_j$,
$1\leq j\leq s$, $\Gamma_1,\dots,\Gamma_s$ are $C^\infty$ 
and in general
position. Denote $\widehat B_m\coloneq\Gamma_{1,\dots,m}$ 
the set of
$(\alpha\:p)\in\Bbb R\Bbb P_n$ such that $l_{\alpha p}$ is 
tangent to $B_m$.
The set $\widehat B_m\subset\Bbb R\Bbb P_n$ may not be a 
hypersurface (see
Theorem 3); however, as Theorem $1'$ claims, it is indeed 
a smooth hypersurface
outside a set of $(n-1)$-dimensional Lebesgue's measure 
zero.
\bh {\rm 3}\endbh
Our second result gives the relation between the 
discontinuity surfaces for
$R(p,\alpha;f)$ and those for $f(x)$; namely, the {\it 
function} $g(x')$ {\it
in the local equation of $\Gamma$}, $x_n=g(x')$, {\it is 
the Legendre transform
of the function $h(\beta)$ which gives the equation of 
$Q_f$}, $q=h(\beta)$.
\par Assume that $q=h(\beta)$, $\beta\in 
U_{\overline\beta}$, where
$U_{\overline\beta}$ is a neighborhood of a point 
$\overline\beta$, $\overline
q=h(\overline\beta)$, and 
$\operatorname{det}h_{ij}(\beta)\not=0$ in
$U_{\overline\beta}$, where
$h_{ij}\coloneq\partial^2h/\partial\beta_i\partial%
\beta_j$. Let $\overline
x'=\nabla h(\overline\beta)$.
\thm{Theorem 2} If $h(\beta)\in C^l(U_{\overline\beta})$, 
$l\geq 2$, then
$Lh=g(x')$, and $g(x')\in C^l(U_{\overline x'})$.
\ethm
\par This result allows one to recover the surfaces of 
discontinuity of $f(x)$
given the surfaces of discontinuity of $R(p,\alpha;f)$.
\bh {\rm 4}\endbh Examples show that the Legendre 
transform $h(\beta)=Lg$ of a
function $g(x')$, $x'\in\Bbb R^{n-1}$, may have domain of 
definition of
dimension less than $n-1$. Since $Q_f$ is a union of 
several varieties of
codimension one in $\Bbb R\Bbb P_n$ (called components 
below), the question
arises: which of the components of $Q_f$ and which of 
their intersections
provide, after applying the generalized Legendre transform 
defined in [RZ1],
parts of $\Gamma=\partial D$ which have codimension one in 
$\Bbb R^n$. The
answer is given in Theorem 3. This theorem describes $Q_f$ 
in terms of
differential geometry of $\Gamma$. Recall that the 
principal curvatures of a
hypersurface $S\subset\Bbb R^n$, which is the graph of a 
function $x_n=g(x')$,
are the eigenvalues of the matrix 
$(g_{ij})\cdot(\delta_{ij}+g_ig_j)^{-1}\cdot
(1+\sum^{n-1}_{i=1}g^2_i)^{-1/2}$, $g_i=\partial 
g/\partial x_i$. One can prove
that {\it if $k$, $k\geq 1$, principal curvatures of a 
hypersurface $S$ vanish
identically, then for every point $P\in S$ there exists an 
affine
$k$-dimensional space $L_P$ such that $P\in L_P\subset S$}.
\thm{Theorem 3} {\rm(a)} Assume that $B_m$ is nonempty. 
Then $m$ principal
curvatures of $\widehat B_m$ vanish identically\RM;
\par{\rm(b)} If $k$ principal curvatures of $\Gamma_1$ 
vanish identically, then
$\widehat\Gamma_1$ has codimension $k+1$ in $\Bbb R\Bbb 
P_n$.
\par Every point of $\widehat\Gamma_1$ is a vertex of a 
cone $K$, which belongs
to $\Gamma_{1j}$, where 
$\Gamma_1\cap\Gamma_j\not=\varnothing$. The directrix
of $K$ is $(k-1)$-dimensional, and this directrix can be 
described as
follows\RM: Take an arbitrary point $P\in\Gamma_1$, and let
$L_k(P)\subseteq\Gamma_1$ be a $k$-dimensional affine 
space containing $P$,
which exists since $k$ principal curvatures of $\Gamma_1$ 
vanish identically.
Let $d_P\coloneq\{(\alpha\:p)\:l_{\alpha p}$ be tangent to 
$\Gamma_j$ at the
points of $L_K(P)\cap\Gamma_j\}$, and let 
$l_{\alpha_0p_0}$ be tangent to
$\Gamma_1$ at the point $P$. The vertex of $K$ is the 
point $(\alpha_0\:p_0)$.
The directrix of $K$ is the set $d_P$.
\ethm
\par The set $Q_f$ is a union of the sets 
$\widehat\Gamma_{i_1\cdots i_k}$,
$Q_f=\bigcup\widehat\Gamma_{i_1\cdots i_k}$ where the 
union is taken over all
combinations of indices $1\leq i_k\leq s$. Theorem 3 gives 
a recipe to select
the components of $Q_f$ which yield after the Legendre 
transform the components
of $\Gamma$ of codimension 1, i.e., hypersurfaces 
$\Gamma_j$ which are parts of
$\Gamma$, $\Gamma=\bigcup^s_{j=1}\Gamma_j$. Note that if a 
component of $Q_f$
has some principal curvatures vanishing identically, then 
its preimage in $\Bbb
R^n$ has codimension greater than one. Therefore, if one 
wishes to recover
hypersurface-type components of $\Gamma$, then one should 
apply the Legendre
transform to those components of $Q_f$, which do not have 
principal curvatures
which vanish identically. Those hypersurfaces $\Gamma_j$ 
which have identically 
vanishing principal curvatures are reconstructed by 
applying the generalized
Legendre transform, which was introduced in [RZ1], to 
high-codimension parts of
$Q_f$ described in Theorem 3(b). The generalized Legendre 
transform was applied
in [Z] to the study of dual varieties in algebraic geometry.
\par It is well known that the Radon transform may be 
considered as a Fourier
integral operator, so it makes sense to study its action 
on the wave front set
of $f$. In [RZ1] we study a relation of the wave front of 
$f$ and the set
$Q_f$.
\heading III. Proofs of Theorems 1 and 2\endheading
\par We sketch the proofs in the simplest case $m=1$, 
$n=2$, but the ideas are
similar in the general case.
\par First we prove that if $D\subset C^\infty$ and $f\in 
C^\infty$, then
$R(p,\alpha;f)\in C^\infty$ on the set $V_f\coloneq\Bbb 
R\Bbb P_n\backslash
Q_f$. Thus, the singularities of $f$ are in the set $Q_f$. 
Second, we prove
that, generically, $Q_f$ is a $C^\infty$ hypersurface in 
$\Bbb R\Bbb P_n$ and
find the equation of this hypersurface.
\par Third, we prove that there exists a neighborhood $U$ 
of a generic point
$(\overline\alpha\:\overline p)$ and an equation 
$\zeta(\alpha\:p)=0$,
$\nabla\zeta\not=0$ in $U$, such that (2) holds.
\par (a) Let us start with the second claim and prove also 
Theorem 2 for $n\geq
2$. Let $\alpha\cdot x-p=0$ be a tangent plane $l_{\alpha 
p}$ to $\Gamma$ at a
point $\overline x\in\Gamma$. Assume that 
$\alpha_n\not=0$, and write
$x_n=\beta\cdot x'-q$, $\beta_i\coloneq-\alpha_i/\alpha_n$,
$q\coloneq-p/\alpha_n$, $x'=(x_1,\dots,x_{n-1})$. Let 
$x_n=g(x')$ be the
equation of $\Gamma$ in a neighborhood $\overline U$ of 
$\overline x$, and
$\operatorname{det}g_{ij}(\overline x')\not=0$. Then 
$\nabla g(x')=\beta$,
$q=\beta\cdot x'-g(x')$. Thus $q=h(\beta)\coloneq Lg$. The 
equation
$q=h(\beta)$ is the equation of $Q_f$ in the inhomogeneous 
coordinates
$(\beta,q)$. One can prove that if $q\in C^s(\overline 
U)$, $s\geq 2$, and
$\operatorname{det}g_{ij}(\overline x')\not=0$, then $h\in 
C^s(U)$, where $U$
is a neighborhood of the point $(\overline\beta,\overline 
q)$, $\nabla
g(\overline x')=\overline\beta$, $\overline 
q=\overline\beta\cdot\overline
x'-g(\overline x')$. Since $L$ is involutive, $g=Lh$. 
Theorem 2 is proved.
\par (b) Let us prove the first claim for $n=2$. Assume 
that $(\alpha\:p)\in
V_f$, i.e., $l_{\alpha p}$ is not tangential to $\Gamma$. 
Write $R(p,\alpha;f)$
as
$$J\coloneq\int^{a_2(q,\beta)}_{a_1(q,\beta)}f(x_1,\beta 
x_1-q)\,dx_1,$$
where $a_i\coloneq a_i(q,\beta)$ are the points of 
intersection of $l_{\alpha
p}$ with $\Gamma$. The integral $J$ is a sum of the 
integrals over the
intervals $(a_1,b),(b,c),(c,a_2)$, where $a_1<b<c<a_2$ and 
$b,c$ do not depend
on $q,\beta$. Obviously the integral over $(b,c)$ is a 
$C^l$ function of
$\beta$ and $q$ if $f\in C^l$, $l\geq 0$. The integrals 
over $(a_1,b)$ and
$(c,a_2)$ are treated similarly.
\par Let us prove that the integral over $(a_1,b)$ is 
$C^l$ function of
$q,\beta$ if $\Gamma,f\in C^l$, $l\geq 2$, and $l_{\alpha 
p}$ is transversal to
$\Gamma$, that is, $\beta\not=g'(a_1)$. It is sufficient 
to prove that
$a_1(q,\beta)\in C^l$. The function $a_1(q,\beta)$ is the 
root of the equation
$q=\beta a_1-g(a_1)$. By the transversality condition 
$\beta-g'(a_1)\not=0$.
Thus, the implicit function theorem implies that the root 
$a_1(q,\beta)\in C^l$
if $g\in C^l$. The first claim is proved.
\par (c) Let us prove the last claim. Let 
$(\overline\alpha\:\overline p)\in
Q_f$ and $(\overline\beta,\overline q)$ be the 
corresponding nonhomogeneous
coordinates. For a generic $(\overline\alpha\:\overline 
p)$ the condition
$g''(\overline x_1)\not=0$ follows from the equation 
$g'(\overline
x_1)=\overline\beta$ and Sard's theorem. We can assume 
therefore that
$g''(\overline x_1)\not=0$. Consequently, the point 
$\overline x_1$ is a
Morse-type (nondegenerate) critical point of the function
$z\coloneq\overline\alpha\cdot x-\overline p$ on 
$\Gamma\cap\overline U$, i.e.,
of the function $-\overline\beta x_1+g(x_1)+\overline q$. 
The part of integral
(1) taken over the complement to $\overline U$ is a 
$C^\infty$-function of
$(\alpha\:p)$ according to $(b)$. It gives $r_2$ in 
formula (2). By the Morse
lemma, there are coordinates $u_1,u_2$ such that the 
equation of $\Gamma$ in
these coordinates is $u_1=0$, the region $D\cap\overline 
U$ is described by the
inequality $u_1\geq 0$, and $z=u_1+u^2_2$ in $\overline 
U$. To study the
singularity of $R(p,\alpha;f)$, take a curve $\gamma$ 
which intersects $Q_f$
transversally, for instance, 
$\gamma=\{(\alpha\:p)\:\alpha=\overline\alpha\}$.
Parameter $p$ gives the position of a point on $\gamma$. 
On $l_{\overline\alpha
p}$ one has $\overline\alpha\boldcdot x-p=0$ and 
$z=\overline\alpha\cdot
x-\overline p$, so $z=p-\overline p$ on 
$l_{\overline\alpha p}$. Thus, $z$ can
be used as a parameter which determines the position of a 
point on $\gamma$;
therefore, the domain of integration in (1) can be 
described by the inequality
$u_1\geq 0$ and the equation $z-u_1-u^2_2=0$. Thus, 
$z-u^2_2=u_1$, so
$-z^{1/2}_+\leq u_2\leq z^{1/2}_+$ since $z=z_+\geq 0$ in 
the integration
region. We have
$$R(p,\overline\alpha;f)=\int_{l_{\overline\alpha
p}}f(x)\mu(dx)=\int_lf_1(u_1,u_2)\mu_1(du)=\int^{z^{1/2}_+
}_{-z^{1/2}_+}f_2(u_2,z_+)\,du_2,$$
where $f_2(u_2,z)$ is a $C^\infty$-function, $l$ is the 
curve given by the
equation $z-u_1-u^2_2=0$ and $u_1\geq 0$, $\mu_1(du)$ 
comes from $\mu(dx)$ via
the Morse lemma change of variables, and the last integral 
comes after an
elimination of $u_1$. From this formula one derives (2). 
Indeed, write
$f_2(u_2,z)$ as a sum of even $f_{\roman e}$ and odd 
$f_{\roman o}$ functions
of $u_2$, $f_e(u_2,z)\in C^\infty$, $f_o(u_2,z)\in 
C^\infty$. Then the integral
$$\int^{z^{1/2}_+}_{-z^{1/2}_+}f_{\roman e}(u_2,z_+
)\,du_2=z^{1/2}_+r_1\quad\text{and}\quad\int
^{z^{1/2}_+}_{-z^{1/2}_+}f_o(u_2,z_+)\,du_2=0,$$
where $r_1\in C^\infty$. The function $r_2$ in formula (2) 
vanishes if $\Gamma$
is strictly convex so that $l_{\overline\alpha p}$ 
intersects $\Gamma$ at
two points only.
\heading IV. Examples\endheading
\par 1. Let $f(x)=1$, $|x|\leq a$, $f(x)=0$, $|x|>a$, 
$x\in\Bbb R^n$, $n\geq
2$, $\hat f(p,\alpha_0)=2\sqrt{a^2-p^2}$, 
$\alpha^0=\alpha|\alpha|^{-1}$. Thus
$p^2/|\alpha|^2=a^2$ is the equation of $Q_f$. In 
$(\beta,q)$ coordinates the
equation of $Q_f$ is $q=\pm a\sqrt{1+\beta^2}$, 
$\beta\in\Bbb R^{n-1}$. Thus
$h(\beta)=\pm a\sqrt{1+\beta^2}$. By Theorem 2 the 
equation $x_n=g(x')$ of the
surface of discontinuity of $f(x)$ is given by 
$g(x')=Lh=\mp\sqrt{a^2-x\prime^2}$.
The equation $x_n=\pm\sqrt{a^2-x^{\prime 2}}$ defines the 
sphere $|x|=a$.
\par 2. Let $f(x)=1$, $b\leq|x|\leq a$, $f(x)=0$, $|x|<b$ 
or $|x|>a$, $0<b<a$,
$n\geq 2$. Then $\hat f(p,\alpha^0)=2\sqrt{a^2-p^2}$, 
$b\leq p\leq a$; $\hat
f(p,\alpha^0)=2(\sqrt{a^2-p^2}-\sqrt{b^2-p^2})$, $0\leq 
p\leq b$; $\hat
f(p,\alpha^0)=0$, $p>a$, $|\alpha^0|=1$. Thus 
$p^2=|\alpha|^2a^2$ and
$p^2=|\alpha|^2b^2$ are the equations of $Q_f$. Taking 
Legendre's transform
yields the surfaces $|x|=a$ and $|x|=b$ of discontinuity 
of $f(x)$.
\par 3. Consider $f(x)=0$ outside of the region $D$ 
bounded by $\Gamma$, where
$\Gamma$ is the union of the curves $x_2=0$ and 
$x_2=x^2_1-1$, and let
$f(x)=1$, $x\in D$. The $R(p,\alpha;f)$ is a function 
whose support is bounded
by the curves $q=\beta$, $q=-\beta$ from below, 
$q=\tfrac14\beta^2+1$ in the
interval $-2\leq\beta\leq 2$, and $q=\beta$, $q=-\beta$ 
for $|\beta|\geq 2$
from above. One can check that on the lines $q=\pm\beta$,
$-\infty<\beta<\infty$, the function $R(p,\alpha;f)$ has a 
singularity of the
type $|z|$ and on the parabola $q=\tfrac14\beta^2+1$ it 
has the singularity of
the type $z^{1/2}_+$. Applying Legendre's transform first 
to the function
$q=\tfrac14\beta^2+1$, $-2\leq\beta\leq 2$, yields the 
parabola $x_2=x^2_1-1$,
$-1\leq x_1\leq 1$; and secondly, applying it to the 
functions $q=\pm\beta$
yields two points $x_1=\pm 1$, $x_2=0$. By Theorem 3, the 
straignt line joining
these two points also belongs to $\Gamma$. Thus $\Gamma$ 
is recovered.
\heading V. Numerical aspects\endheading
\par The RT of $f(x)$ is usually given with an error. 
Hence, the first
numerical problem is to calculate the function $h(\beta)$ 
which gives the
equation of the set $Q_f$ of the singularities of RT given 
the noisy
measurements of the RT. The second numerical problem is to 
calculate
$Lh=g(x')$. Calculation of the Legendre transform of a 
function $h(\beta)$
known with errors is a well-posed problem, at least in the 
case when
$\operatorname{det}g_{ij}(x')\not=0$. It is proved in 
[RSZ] that if a function
$g_\delta(x')$ given such that 
$|g_\delta(x')-g(x')|<\delta$, $g_\delta(x')$ is
not necessarily in $C^2$ but is continuous, then one can 
calculate $Lg$ with
the accuracy $O(\delta)$ as $\delta\to0$. This means that 
a stable method is
given in [RSZ] for calculating the Legendre transform of 
noisy data. See also
[R5]. Our result in part 3 of \S II has an interesting 
connection with the
envelopes theory [T, Zl].
\heading Acknowledgments\endheading A. G. Ramm thanks ONR, 
NSF, and USIEF for
support. The research of A. I. Zaslavsky was supported in 
part by a grant from
the Ministry of Science and the ``Ma-agara''-special 
project for absorption of
new immigrants, in the Department of Mathematics, Technion.

\enddocument